\documentclass[]{article}

\usepackage[margin=1.5in]{geometry}
\usepackage{latexsym}
\usepackage{amsmath}
\usepackage{amssymb,amsfonts,amscd}
\usepackage{graphicx}
\usepackage{epstopdf}
\graphicspath{{./diagrams/}}
\usepackage{multirow}
\usepackage{hyperref}
\usepackage{bm,color,mathdots}

\newtheorem{theorem}{Theorem}

\newtheorem{ex}[theorem]{Example}

\usepackage{algorithm}
\usepackage[noend]{algorithmic}
\floatname{algorithm}{Procedure}

\usepackage{marginnote}

\newcommand{\MYcomment}[1]{\textcolor{red}{$\star$}\marginnote{#1}}

\makeatletter
\let\oldmarginnote\marginnote
\renewcommand*{\marginnote}[1]{%
   \begingroup%
   \ifodd\value{page}
     \if@firstcolumn\reversemarginpar\fi
   \else
     \if@firstcolumn\else\reversemarginpar\fi
   \fi
   \oldmarginnote{#1}%
   \endgroup%
}
\makeatother

\newcommand{\bC}{{\mathbb C}}

\newcommand{\bP}{{\mathbb P}}
\newcommand{\bQ}{{\mathbb Q}}
\newcommand{\bR}{{\mathbb R}}

\newcommand{\bZ}{{\mathbb Z}}

\newcommand\sA{{\mathcal A}}
\newcommand\sB{{\mathcal B}}

\newcommand\sD{{\mathcal D}}

\newcommand\sI{{\mathcal I}}

\newcommand\sL{{\mathcal L}}

\newcommand\sV{{\mathcal V}}

\newcommand\Var{\sV}

\def\rank{{\mathop{\rm rank~}\nolimits}}
\def\Null{{\mathop{\rm null~}\nolimits}}
\def\dnull{{\mathop{\rm dnull}\nolimits}}
\def\Iso{{\mathop{\rm Iso}\nolimits}}

\newcommand{\suchthat}{\>|\>}
\newcommand{\leftsuchthat}{\>\left|\>}

\newcommand{\rrad}[1]{\sqrt[\bR]{#1}}
\newcommand{\rvar}[1]{\sV_\bR(#1)}
\newcommand{\cbrt}[1]{\sqrt[3]{#1}}

\usepackage{tikz}
\usetikzlibrary{arrows}

\begin{document}
	
	\title{Numerically validating the completeness of the real solution set of a system of polynomial equations}

\author{
Daniel A. Brake and Jonathan D. Hauenstein and Alan C. Liddell
}
	
	
	\maketitle
	
	\begin{abstract}
\noindent 
Computing the real solutions to a system of polynomial equations is a challenging problem, 
particularly verifying that all solutions have been computed. 
We describe an approach
that combines numerical algebraic geometry and sums of squares programming
to test whether a given set is ``complete'' with respect to the real solution
set.  Specifically, we test whether the Zariski closure of that given set is 
indeed equal to the solution set of the real radical of the ideal 
generated by the given polynomials. 
Examples with finitely and infinitely many real solutions are provided, 
along with an example having polynomial inequalities.
	\end{abstract}
	
	
	
	
	\section{Introduction}\label{Sec:Intro}
	
Numerical methods provide approximate solutions to 
continuous problems.  For example, numerical algebraic geometry uses 
numerical methods to compute approximations of the solutions to 
systems of polynomial equations.
Due to the potential for error in numerical approaches, techniques have 
been developed for certifying aspects of numerically computed results.  This article uses
sums of squares programming to validate that 
a {\em complete} real solution set has been computed,
that is, the Zariski closure of the given set is equal to the Zariski
closure of the set of all real solutions.  

A typical situation where one may need to test 
the completeness of a real solution set
is computing critical points.  For example, \S~\ref{Sec:Phi4N3}
considers computing the critical points of a potential energy landscape.
In such situations, local numerical methods, e.g., \cite{DoyeWales2002,mehta2014communication},
exist for locating real critical points. Our approach provides a 
global stopping criterion for validating that all real solutions have been
identified.

A related situation is the computation of the real critical points of a
projection of a solution set used in the numerical decomposition of 
real curves and surfaces \cite{RealSurfaces2,RealSurfaces,brake2014bertini_real,RealCurves}.
The failure to correctly compute the set of real solutions leads to a 
failure in the decomposition of the real component.
Hence, correct and complete computation of sets of real solutions 
is paramount to correctly computing the decomposition.

One approach for certifying the existence of real solutions 
is based on the local analysis of Newton's method using 
Smale's $\alpha$-theory \cite{smale1986newton} developed
in \cite{Cite:alphaCert}.  Building on \mbox{$\alpha$-theory},
there are methods for certifying smooth continuous paths 
for Newton homotopies \cite{Cite:cadenza,hauenstein2014certified} 
and general homotopies \cite{beltran2013robust}.
For example, if a smooth path is defined by a real system of equations
which has a real starting point, then the endpoint of the path must also be real.
	
From an algebraic viewpoint, the radical of an ideal generated by a
given collection of polynomials consists of all polynomials that vanish on the solution
set of the given polynomials.  There are several algorithms for computing the 
radical of a zero-dimensional ideal -- some numerical, 
e.g.,~\cite{janovitz2012computation,lasserre2009prolongation,lasserre2009unified} 
and some symbolic, e.g., \cite{becker1996radical,gianni1988grobner}.
When there are infinitely many solutions, one can reduce to the zero-dimensional
case, for example, via \cite{gianni1988grobner,krick1991algorithm}.
	
The {\em real} radical of an ideal generated by a given collection of 
polynomials with real coefficients consists of all polynomials 
that vanish on the {\em real} solution set of the given polynomials.
There have been several proposed methods for computing the real radical of an ideal.
Some are symbolic, e.g.,~\cite{becker1993computation} based on the primary decomposition (see also \cite{neuhaus1998computation,spang2007computation,xia2002algorithm,zeng1999computation}).
Others are numerical, based on moment matrices when the number of real solutions 
is finite, e.g., \cite{lasserre2013moment,LLR08,lasserre2009prolongation,lasserre2009unified, laurent2012approach}.  A promising approach for computing the real radical 
when there are infinitely many real solutions was developed in \cite{ma2014certificate}
providing a stopping criterion for verifying that a Pommaret basis has been computed.
Other methods for computing real solutions include 
computing a point on each semi-algebraically connected component of 
the real solution set, e.g., \cite{aubry2002real,bank2001polar,hauenstein2013numerically,safey2003polar},

As discussed in \cite{ma2014certificate}, one key issue related to 
computing the real radical using semidefinite programming with moment matrices
is knowing when the generated polynomials form a basis for the real radical.  
In our approach, we first compute a set $S$ which is a subset of 
the Zariski closure of the real solution set.
Then, we compute polynomials that vanish on~$S$.  
Finally, for each of the computed polynomials, 
we use sums of squares programming to verify that it is indeed in 
the real radical.  Since the polynomials can be validated independently,
we can easily parallelize this part of the computation.
Since $S$ is contained in the Zariski closure of the real solution set,
every polynomial contained in the real radical vanishes on $S$.
Conversely, if every polynomial that vanishes on $S$ is contained 
in the real radical, we know that a generating set for the real radical
has been computed.  Hence, $S$ is complete since the Zariski closure of 
$S$ is equal to the Zariski closure of the real solution set of 
the original system of equations, i.e., the
solution set of the real radical.

We perform these computations numerically.
From the numerical output, one could then aim to produce 
exact representations of the polynomials, e.g., via~\cite{ExactnessRecovery}.
This would typically require field extensions, which is one of the pitfalls 
of using purely symbolic methods to compute real radicals.
As an illustrative example, 
consider the polynomial $f(x) = x^3 - 2$
having rational coefficients, i.e., $f\in\bQ[x]$.
Since $f = 0$ has one real solution, namely $x = \cbrt{2}$,
the real radical of the ideal generated by $f$ is 
$\langle x - \cbrt{2}\rangle$ 
which is generated by a polynomial not in $\bQ[x]$.
From a numerical approximation of $\cbrt{2}$,
exactness recovery methods, e.g., \cite{ExactnessRecovery},
allow one to determine that exact results could be obtained by 
working over the coefficient field $\bQ[\cbrt{2}]$.  

The remainder of the article is as follows.
Section~\ref{Sec:Ideals} focuses on 
radicals, irreducible decomposition, and Zariski closures.  
Real radicals, sums of squares, and semidefinite programming
are discussed in Section~\ref{Sec:SemiSoS}.
Section~\ref{Sec:Sampling} considers generating a subset $S$
contained in the Zariski closure of the real solution set, 
including a discussion on finding and sampling positive-dimensional components.
From the set~$S$, interpolation is used to compute a candidate set of generators for 
the real radical as described in Section \ref{Sec:Interpolate}.
Section~\ref{Sec:Certification} presents a 
criterion for showing that a set $S$ is complete
with respect to the real radical.
Section~\ref{Sec:Inequalities} considers 
the real solution set for collection of equations and inequalities. 
Several examples are presented in Section~\ref{Sec:Examples} 
and we conclude in Section~\ref{Sec:Conclusion}.

\subsection*{Acknowledgments}

The authors thank Mohab Safey El Din, Charles Wampler,
and Lihong Zhi for discussions 
related to real solution sets.

	
	\section{Zariski closure and radicals}\label{Sec:Ideals}

Let $f_1,\dots,f_k\in\bC[x_1,\dots,x_n]$ and consider
the ideal generated by these polynomials, namely $I = \langle f_1,\dots,f_k\rangle$.
The polynomials $f = \{f_1,\dots,f_k\}$ and the corresponding ideal 
$I = \langle f\rangle$ define
the same solution set in $\bC^n$, namely
$$\Var_\bC(f) = \Var_\bC(I) = \{x\in\bC^n\suchthat f_i(x) = 0 \mbox{~for~} i = 1,\dots,k\}$$
A set $A\subset\bC^n$ is called an {\em algebraic set} if
there is a collection of polynomials $g\subset\bC[x_1,\dots,x_n]$ such
that $A = \Var_\bC(g)$.  The algebraic set $A$ is {\em irreducible}
if there does not exist algebraic sets $A_1,A_2\subsetneq A$ with 
$A = A_1\cup A_2$.  Given an algebraic set $A$, there exists a unique collection (up to relabeling) of irreducible algebraic sets $X_1,\dots,X_\ell$ such that
$$A = \bigcup_{i=1}^\ell X_i \hbox{~~and~~} X_j\not\subset \bigcup_{i\neq j} X_i.$$
Each $X_i$ is called an {\em irreducible component} of $A$.  

In numerical algebraic geometry, an irreducible algebraic set is
represented by a {\em witness set}, see, e.g., \cite[Chap.~13]{SW05}.
A {\em numerical irreducible decomposition} for an algebraic set $A$ 
is a collection of witness sets for the irreducible components of $A$.
Such a decomposition can be computed using various algorithms, 
e.g., \cite{ldt,RegenCascade,Cascade,NumAlgGeom}.

For any subset $T\subset\bC^n$, the ideal generated by $T$ is
$$I(T) = \{f\in\bC[x_1,\dots,x_n]~|~f(t) = 0 \hbox{~for all~}t\in T\}.$$
The {\em Zariski closure} of $T$ is the algebraic set $\overline{T} = \Var_\bC(I(T))$,
which is the intersection of all algebraic sets that contain~$T$.

For an ideal $I$, the {\em radical} of $I$ is $\sqrt{I} = I(\Var_\bC(I))$
which can be described algebraically as
$$\sqrt{I} = \{p\in\bC[x_1,\dots,x_n]~|~p^\alpha \in I \hbox{~for some~}\alpha\in\bZ_{> 0}\}.$$
	
	
	\section{Real radical \& sums of squares}\label{Sec:SemiSoS}

Many of the topics 
from \S~\ref{Sec:Ideals} have analogous statements over $\bR$.
Let $f_1,\dots,f_k\in\bR[x_1,\dots,x_n]$ with $f = \{f_1,\dots,f_k\}$
and $I = \langle f\rangle$.  The set of solutions in $\bR^n$ is
$$\mbox{\scriptsize $\Var_\bR(f) = \Var_\bR(I) = \{x\in\bR^n\suchthat f_i(x) = 0 \mbox{~for~} i = 1,\dots,k\} = \Var_\bC(I)\cap\bR^n.$}$$
The {\em real radical} of~$I$ is $\sqrt[\bR]{I} = I(\Var_\bR(I))$
which can also be described algebraically as
\begin{equation}\label{eq:RealRadical}
\sqrt[\bR]{I} = \left\{p\in\bR[x]\leftsuchthat\begin{array}{l}
p^{2\alpha} + \sum_{j=1}^\ell g_j^2 \in I  \\
\mbox{~~~for some~}\alpha \in \bZ_{> 0}, g_j\in\bR[x] \end{array} \right\}\right..
\end{equation}

\begin{ex}\label{Ex:Illustrative}
For $f(x) = x^3 - 2$ and $I = \langle f\rangle$, we have:
\begin{itemize}
\item $\Var_\bC(I) = \{\cbrt{2},\omega \cbrt{2}, \omega^2 \cbrt{2}\}$ and $\Var_\bR(I) = \{\cbrt{2}\}$,
\item $\sqrt{I} = I$, and
\item $\sqrt[\bR]{I} = \langle x - \cbrt{2} \rangle$
\end{itemize}
where $\omega$ is the primitive cube root of unity.  
In particular,
$$
(x - \cbrt{2})^4 + (\sqrt{3}x^2 - \sqrt{3} \cbrt{4})^2 = 4(x^3 - 2)(x - \cbrt{2}) \in I.
$$
\end{ex}

The algebraic description of the real radical $\sqrt[\bR]{I}$ 
presented in \eqref{eq:RealRadical} shows that this
definition depends on {\em sums of squares}.
A polynomial $s\in\bR[x_1,\dots,x_k]$ is called a sum of squares 
if $s = \sum_{j=1}^\ell g_j^2$ for some $g_1,\ldots,g_\ell\in\bR[x_1,\dots,x_k]$.
Clearly, every polynomial that is a sum of squares
has even degree.  

The polynomials of even degree that are sums of squares
are characterized by {\em positive semidefinite} matrices.
A symmetric matrix $M \in \bR^{m \times m}$ is positive semidefinite 
if, for all $y\in\bR^m$, \mbox{$y^T M y \geq 0$}.
This condition is equivalent to all eigenvalues of $M$ being nonnegative.
We will write $M\succeq 0$ if $M$ is positive semidefinite.

Let $s\in\bR[x_1,\dots,x_k]$ be a polynomial of degree $2d$ 
and~$X_d$ be 
the vector of all monomials in $x_1,\dots,x_n$ of degree at most~$d$.
Hence, there exists a symmetric matrix $C$ such that
\begin{equation}\label{eq:C}
s(x) = X_d^T\cdot C\cdot X_d.
\end{equation}
The polynomial $s$ is a sum of squares if and only if there is a positive
semidefinite matrix $C$ such that \eqref{eq:C} holds.

\begin{ex}\label{Ex:Illustrative2}
As shown in Ex.~\ref{Ex:Illustrative}, the quartic polynomial 
$s(x) = 4(x^3-2)(x-\cbrt{2})$ is a sum of squares.  Let
$$X_2 = \left[\begin{array}{c} 1 \\ x \\ x^2 \end{array}\right] 
\mbox{~and~} C = \left[\begin{array}{ccc}
8\cbrt{2} & -4 & -2\cbrt{4} \\ 
-4 & 4\cbrt{4} & -2\cbrt{2} \\
-2\cbrt{4} & -2\cbrt{2} & 4
\end{array}\right].$$
It is easy to verify that $C\succeq 0$ and $s(x) = X_2^T\cdot C\cdot X_2$.
\end{ex}

For a given polynomial $s$ of degree $2d$, the set 
of symmetric matrices $C$ such that \eqref{eq:C}
holds is a linear space.  Hence,
testing that a polynomial is a sum of squares 
can be accomplished by solving a semidefinite feasibility problem.

\begin{ex}\label{Ex:Illustrative3}
Continuing with $s(x) = 4(x^3-2)(x-\cbrt{2})$ from Ex.~\ref{Ex:Illustrative2}, consider
the linear space
$$\sL = \left\{\left[\begin{array}{ccc}
s_{00} & s_{01} & s_{02} \\ s_{01} & s_{11} & s_{12} \\ s_{02} & s_{12} & s_{22} 
\end{array}\right] \leftsuchthat 
{\scriptsize
\begin{array}{rcl} 
s_{00} & = & 8\cbrt{2} \\
2s_{01} & = & -8 \\
2s_{02} + s_{11} & = & 0 \\
2s_{12} & = & -4\cbrt{2} \\
s_{22} & = & 4 \end{array}}
\right\}\right..$$
Since $s(x) = X_2^T \cdot C\cdot X_2$ if and only if $C\in\sL$,
it follows that~$s$ is a sum of squares if and only if there exists
$C\in\sL$ such that $C\succeq0$, which is 
a semidefinite feasibility problem.
\end{ex}

Since the task of converting between sums of squares problems
and semidefinite programming problems can be arduous,
we utilize the software package {\tt SOSTOOLS} \cite{Cite:SOSTools}.

Given a polynomial $p\in\bR[x_1,\dots,x_n]$, we can 
decide if \mbox{$p\in\sqrt[\bR]{I}$} using \eqref{eq:RealRadical}.
That is, 
$p\in\sqrt[\bR]{I}$ if and only if there exists $\alpha\in\bZ_{> 0}$
and $h_1,\dots,h_k,g_1,\dots,g_\ell\in\bR[x_1,\dots,x_n]$ such that
$$p^{2\alpha} + \sum_{j=1}^\ell g_j^2 = \sum_{i=1}^k h_i f_i$$
which is equivalent to requiring that 
\begin{equation}\label{eq:SOS}
-p^{2\alpha} + \sum_{i=1}^k h_i f_i\mbox{~is a sum of squares}.
\end{equation}
Thus, given a polynomial $p\in\bR[x_1,\dots,x_n]$, one can test
if $p\in\sqrt[\bR]{I}$ by solving a 
semidefinite feasibility problem.  The construction of such
polynomials $p$ used for testing 
is based on computing points in $\Var_\bR(I)$,
which is discussed next.
	
	
	\section{Generating a candidate set}\label{Sec:Sampling}

The key aspect of our approach is to first produce a superset of
the real radical ideal.  This is accomplished by
computing a set $S\subset\overline{\Var_\bR(I)}$.
In particular, if $S\subset\overline{\Var_\bR(I)}$, then 
$\sqrt[\bR]{I}\subset I(S)$.
We then aim to show that $I(S) = \sqrt[\bR]{I}$.
Since our approach is dependent on the ability to 
generate~$S$, we discuss several possible methods for procuring $S$.
	
\subsection{Approaches for locating real solutions}	
	
A classical approach for attempting to find a real solution is 
to use Newton's method or related variants, see, e.g., \cite{kelley2003solving}.
For a polynomial system with real coefficients, 
if the initial point is real, then every solution obtained from 
Newton's method is also real.  Of course, there are many challenges 
associated with finding real solutions using Newton's method, 
particularly when $\Var_\bC(f)$ is not a complete intersection
or the real solutions are singular with respect to $f$.  
That is, problems can occur with Newton's method,
e.g., divergence, if the dimension of the solution set is less than 
dimension of the null space of the Jacobian at the solution~\cite{GO2,GO1}.  
Nonetheless, heuristic techniques such as 
damping methods, reusing Jacobians for several iterations, 
or using chord or secant methods can be utilized \cite{kelley2003solving}.

Another approach for computing real solutions is to utilize
numerical optimization techniques.  Standard iterative 
techniques include those based on nonlinear least squares approaches
such as the Levenberg-Marquardt algorithm and alternating least squares \cite{kelleyOptimization}.  Other standard methods in optimization
include the worker bees method, genetic algorithms, and the Nelder-Mead 
method, see, e.g., \cite{chong2013introduction}.

Critical point methods combine optimization and polynomial
system solving techniques.  For example, Seidenberg~\cite{S54}
considered the critical points of the distance function between
the set of real solutions and a given real point~$y^*$ that was not a solution.  
The set of all such critical points contains a point on every connected
component of the real solution set \cite{ARS02,RRS00,S54}.
By utilizing homotopy continuation, one can compute a finite subset
of critical points containing a point on every connected component \cite{hauenstein2013numerically}.  Moreover, one can then sample
more real points by moving $y^*$.  

Rather than compute all critical points, one can attempt to compute
the closest critical point to the given $y^*$.  This
can be accomplished using a classical optimization approach
such as the gradient descent method or a homotopy-based
approach called gradient descent homotopies \cite{GH15}.
By testing at many values of $y^*$, one aims to quickly generate
many real solutions, e.g., as shown in \cite[Fig.~3]{GH15}.

Other so-called ``local'' solving methods exist for finding
real solutions, which have been used in various disciplines.
Some examples include techniques in theoretical chemistry,
e.g., \cite{DoyeWales2002,mehta2014communication,mehta2015exploring} 
and solving power-flow equations in electrical engineering,
e.g., \cite{LW,MaThorp}.  

\subsection{Real solutions and isosingular sets}\label{Ssec:Isosingular}

After a real solution has been located, 
one can now try to extract additional information about the 
geometry of the solution set near this point.  One approach
is to compute a local irreducible decomposition using
local witness sets~\cite{BHS} to see if local structure
provides insight into the components of the real solution set 
passing through the computed real point.  Another approach is to utilize
{\em isosingular sets}~\cite{HW13}, which may also
help in improving the numerical stability of interpolation,
described in the next section.

Let $f_1,\dots,f_k$ be polynomials and $z\in\Var_\bC(f)$.
Let $Jf(z)$ be the Jacobian matrix of $f$ evaluated at $z$.
For an integer $\ell$, let $\det_\ell Jf(z)$ be the collection
of all $(\ell+1)\times(\ell+1)$ minors of~$Jf(z)$.
Thus, $\det_\ell Jf(z) = 0$ if and only if $\rank Jf(z) \leq \ell$.
For a polynomial system $g$, let $\dnull(g,z) = \dim \Null Jg(z)$.
The {\em deflation sequence} of $z$ with respect to $f$ is defined by
$$d_i(f,z) = \dnull(\sD^i(f,z),z) \hbox{~for $i\in\bZ_{\geq0}$}$$
where $\sD^0(f,z) = f$ and 
$$\sD^i(f,z) = \left[\begin{array}{c} \sD^{i-1}(f,z) \\ \det_{d_{i-1}(f,z)} J\sD^{i-1}(f,z) \end{array}\right].$$
The deflation sequence is a nonincreasing sequence of nonnegative integers
and thus has a limit, say $d_\infty(f,z)\geq0$, called the 
{\em isosingular local dimension} of $z$ with respect to $f$.

If $X(f,z)$ is the Zariski closure of all points in $\Var_\bC(f)$
which have the same deflation sequence with respect to $f$ as $z$, 
then \cite[Lemma~5.14]{HW13} yields that there is a unique irreducible 
component of $X(f,z)$ which contains $z$, denoted $\Iso_f(z)$,
called the {\em isosingular set} of $z$ with respect to $f$.
In particular, $d_\infty(f,z) = \dim \Iso_f(z)$.

Suppose that $z\in\Var_\bR(f)\subset\bR^n$.  Since $z$ is a smooth
point on the irreducible set $\Iso_f(z)$, we have 
\mbox{$\Iso_f(z)\cap\bR^n\subset\Var_\bR(f)$}
and 
$\Iso_f(z)=\overline{\Iso_f(z)\cap\bR^n}
\subset\overline{\Var_\bR(f)}$.
That is, if \mbox{$I = \langle f\rangle$}, 
$$\Iso_f(z)\subset\Var_\bC(\sqrt[\bR]{I}) \hbox{~~and~~}
\sqrt[\bR]{I}\subset I(\Iso_f(z)).$$
The isosingular local dimension 
is a lower bound on the local real dimension at $z$
which is sharp if $z$ is a smooth point on a unique irreducible
component of $\Var_\bC(\sqrt[\bR]{I})$.  
Moreover, if $d_\infty(f,z) > 0$, 
we can use standard sampling techniques
in numerical algebraic geometry, see, e.g., \cite[\S~8.3]{Cite:bertini},
applied to $\Iso_f(z)$ to produce an arbitrary number
of additional points for which 
polynomials in $\sqrt[\bR]{I}$ must vanish.

Additionally, by using isosingular sets and numerical algebraic geometry,
we can utilize standard membership tests, see, e.g., \cite[\S~8.4]{Cite:bertini},
to determine if a newly found point $x\in\Var_\bR(I)$ is already
contained in the set $S$.

	
	\section{Interpolation} 
	\label{Sec:Interpolate}
	
From the set $S\subset\overline{\Var_\bR(I)}$ 
constructed in \S~\ref{Sec:Sampling}, the next 
task is to compute a collection of polynomials which vanish
on~$S$.  Testing whether $I(S)$ is equal to $\sqrt[\bR]{I}$
is described in \S~\ref{Sec:Certification}.  Here,
we describe computing a basis for $I(S)$ via interpolation.

Suppose that $T\subset\bC^n$ is a finite set 
such that $I(T)$ is generated by
real polynomials and $d\geq1$.  Let $\sB$ form a basis for 
the finite-dimensional vector space of all polynomials in $n$ variables
with real coefficients of degree at most~$d$, namely $\bR[x_1,\dots,x_n]_{\leq d}$.
The linear space of polynomials of degree at most $d$ 
in $I(T)$, denoted $I(T)_{\leq d}$,
is (isomorphic~to) the null space of matrix $M$ where $M_{ij} = \beta_j(t_i)$,
i.e., the evaluation of the $j^{\rm th}$ basis element 
$\beta_j\in\sB$ at the $i^{\rm th}$~point $t_i\in T$.
If $S$ is a finite set, then we simply take $T = S$.  Otherwise,
one can take $T$ to be a finite set consisting of sufficiently
many points on each component described by $S$.  The number of 
sample points needed on each component can be {\em a~priori} bounded
based on the dimension of $\bR[x_1,\dots,x_n]_{\leq d}$.  
One can also algorithmically bound the number of sample points
needed per component simply by continuing to add sample 
points from each component to $T$ until the rank 
of the associated matrix $M$ stabilizes.  

As shown in \cite{GHPS}, one can rescale each row independently
to improve the conditioning of interpolation.  Moreover, for positive-dimensional
components, sampling points that are spread out 
over the component using
numerical algebraic geometry as in \S~\ref{Ssec:Isosingular}
also helps to improve conditioning.  

\begin{ex}\label{example:PlanarSpheroloid}
The solution set of the polynomial system
	\begin{equation}\label{Ex:PlanarSpheroloid}
	f = \{x^2 + y^2 + z^2 - 1,~x^2 + y^2 + z - 1,~x\}
	\end{equation}
consists of the three points
	\[
	\sV_\bC(f) = \rvar{f} = \{(0, 1, 0), (0, -1, 0), (0, 0, 1)\}
	\]
where the point $(0, 0, 1)$ has multiplicity two with respect to~$f$.  

To illustrate, for $d = 2$, we choose the monomial basis
		\[
		\sB = \{1, x, y, z, x^2, xy, xz, y^2, yz, z^2\}
		\]
		for $\bR[x,y,z]_{\leq 2}$ with $S = T = \Var_\bR(f)$ 
		where $M$ is
		\[
		\begin{array}{r|rrrrrrrrrr}
		& 1 & x & y & z & x^2 & xy & xz & y^2 & yz & z^2 \\
		\hline
		(0,1,0) & 1  &  0  &  \phantom{-}1  &  0  &  0  &  0  &  0  &  1  &  0  &  0 \\
		(0,-1,0) & 1  &  0  &            -1  &  0  &  0  &  0  &  0  &  1  &  0  &  0 \\
		(0,0,1) & 1  &  0  &  \phantom{-}0  &  1  &  0  &  0  &  0  &  0  &  0  &  1
		\end{array}.
		\]
A basis for $\Null M$ is given by the columns of the matrix
		\[
\mbox{\small $		
		\left[\begin{array}{rrrrrrr}
 0& 0& 0& 0& -1& 0&  0\\
 1& 0& 0& 0&  0& 0&  0\\
 0& 0& 0& 0&  0& 0&  0\\
 0& 0& 0& 0&  1& 0& -1\\
 0& 1& 0& 0&  0& 0&  0\\
 0& 0& 1& 0&  0& 0&  0\\
 0& 0& 0& 1&  0& 0&  0\\
 0& 0& 0& 0&  1& 0&  0\\
 0& 0& 0& 0&  0& 1&  0\\
 0& 0& 0& 0&  0& 0&  1 
\end{array}\right]$}
		\]
		corresponding to the polynomials
		\begin{align*}
		 x,~x^2,~xy,~xz,~y^2 + z - 1,~yz,~z^2 - z
		\end{align*}
		which form a basis for the linear space $(\sqrt[\bR]{I})_{\leq 2}$.
Note that since each polynomial $f_i$ has degree at most $2$, 
each $f_i$ is contained in the linear span of these polynomials.
		
For illustrative purposes, we selected a monomial basis.
In practice, the choice of basis should be made based on numerical conditioning.
\end{ex}

For $d\gg0$, we know $I(S) = \langle I(S)_{\leq d}\rangle$.
If $S$ is a finite set, then one can determine an upper bound
on $d$ such that $I(S)$ is generated by $I(S)_{\leq d}$.  
In particular, the function 
$$c\mapsto\dim \bR[x_1,\dots,x_n]_{\leq c} - \dim I(S)_{\leq c}$$ 
is the Hilbert function of $I(S)$.  If $r$ is the minimum
such that $|S| = \dim \bR[x_1,\dots,x_n]_{\leq r} - \dim I(S)_{\leq r}$,
i.e., the index of regularity, then one knows that $I(S)$ is either generated by $I(S)_{\leq r}$ or 
$I(S)_{\leq r+1}$. 
In fact, $I(S)_{\leq r}$ generates $I(S)$ if and only if 
$\langle I(S)_{\leq r}\rangle_{\leq r+1} = I(S)_{\leq r+1}$, i.e.,
the Hilbert function of $J = \langle I(S)_{\leq r}\rangle$ in
degree $r+1$ is also equal to $|S|$.  

\begin{ex}
Continuing with Ex.~\ref{example:PlanarSpheroloid}, since
$$\dim\bR[x,y,z]_{\leq 2} - \dim I(S)_{\leq 2} = 10 - 7 = 3 = |S|,$$
one can easily verify that $I(S)$ is generated by $I(S)_{\leq 2}$, i.e.,
$$\sqrt[\bR]{I} = \langle x, y^2 + z - 1, yz, z^2 - z\rangle.$$
\end{ex} 

\begin{ex}
The Hilbert function for the ideal $I(S)$ where $S = \{(0,0),(0,1),(1,0)\}$
is $1,3,3,\dots$ so that $I(S)$ is either generated
by $I(S)_{\leq 1}$ or $I(S)_{\leq 2}$.  Since $I(S)_{\leq 1} = \{0\}$,
we know that $I(S)_{\leq 2}$ must generate $I(S)$.
\end{ex}

When $S$ is infinite, we aim to reduce our computations to 
standard computations performed over $\bC$ as summarized 
in~\S~\ref{Sec:Ideals}.
In particular, by using isosingular sets as discussed 
in~\S~\ref{Ssec:Isosingular}, 
we can actually assume that $S = \overline{S}$ and 
that we have a numerical irreducible decomposition of $S$.
Hence, we simply need to compute $d$ large enough so that
$S$ and $\Var_\bC(I(S)_{\leq d})$ have the same irreducible components 
so that $S = \Var_\bC(I(S)_{\leq d})$.
Hence, $I(S) = \sqrt{\langle I(S)_{\leq d}\rangle}$.  

		
		\section{Validation}\label{Sec:Certification}
		
After computing polynomials which vanish
on $S$, the last step is to verify that they indeed 
lie in the real radical ideal.  Since $S\subset
\overline{\Var_\bR(I)} = \Var_\bC(\sqrt[\bR]{I})$,
we know that $\sqrt[\bR]{I} \subset I(S)$.  
Let $g_1,\dots,g_\ell \in\bR[x_1,\dots,x_n]$ such that
$I(S) = \langle g_1,\dots,g_\ell\rangle$.
If each $g_i \in \sqrt[\bR]{I}$, then we know $I(S) = \sqrt[\bR]{I}$.

Let $I = \langle f_1,\dots,f_k\rangle$.  
For a given $p\in\bR[x_1,\dots,x_n]$, 
we know $p\in\sqrt[\bR]{I}$ if and only if there
exists $\alpha\in\bZ_{>0}$ and 
$h_1,\dots,h_k\in\bR[x_1,\dots,x_n]$ such that \eqref{eq:SOS} holds.
In particular, \eqref{eq:SOS} holds for each $p = g_i$ 
if and only if $I(S) = \sqrt[\bR]{I}$.

If $p\not\in\sqrt[\bR]{I}$, then, for every $\alpha\in\bZ_{>0}$,
\eqref{eq:SOS} does not hold.
Since we can only test finitely many $\alpha$,
an {\em a priori} upper bound on the largest possible value for $\alpha$ 
would be useful for validating that $\sqrt[\bR]{I}\subsetneq I(S)$.
However, without such a bound, we simply keep 
searching for new points in $\Var_\bR(I)$.
If $p\not\in\sqrt[\bR]{I}$, then there must exist a point $x\in\Var_\bR(I)$
such that $p(x)\neq0$.  In fact, there is an irreducible component
$X\subset\overline{\Var_\bR(I)} = \Var_\bC(\sqrt[\bR]{I})$ 
such that $p(x)\neq0$ for every $x$ in a dense open subset of $X$.

With this setup, Procedure~\ref{alg:cert} summarizes our 
complete approach.  If this procedure returns {\sc False}, 
then we either look to add other real solutions to $S$ using \S~\ref{Sec:Sampling} 
or try again with a larger upper bound $\alpha_{\text{max}}$.
We note that, from a practical point-of-view, 
the computations for validation over $\bR$ 
can be simplified by first performing standard computations over $\bC$.  
For example, since $\sqrt[\bR]{I} = \sqrt[\bR]{\sqrt{I}}$, we
could replace $f_1,\dots,f_k$ with a Gr\"obner basis
for $\sqrt{\langle f_1,\dots,f_k\rangle}$.
				
		\begin{algorithm}[h]
			\caption{Validating Real Solution Sets}
			\label{alg:cert}
			\begin{algorithmic}[1]
				\REQUIRE{Polynomials $f = \{f_1,\dots,f_k\} \subset \bR[x_1, \dots, x_n]$ and integer $\alpha_{\text{max}} \in \bZ_{\geq 0}$.}
				\ENSURE{A set $S\subset\overline{\Var_\bR(I)}$ and boolean
				which is {\sc True} if $I(S) = \sqrt[\bR]{I}$ can be validated
				with $\alpha\leq\alpha_{\text{max}}$ 				
				where \mbox{$I = \langle f_1,\dots,f_k\rangle$}, otherwise {\sc False}.}
				\STATE Generate a candidate set $S$ as described in \S~\ref{Sec:Sampling}.
				\STATE Compute polynomials $g_1,\dots,g_\ell\in\bR[x_1,\dots,x_n]$ which generate $I(S)$ as described in \S~\ref{Sec:Interpolate}.\\
				\STATE (Optional) Replace $f$ with a Gr\"obner basis for $\sqrt{I}$.
				\FOR{$m = 1,\dots,\ell$}
				  \STATE Initialize $\alpha:=1$ and $success := $ {\sc False}.
				  \WHILE{$success = $ {\sc False}}
				    \IF{there exists $h_i\in\bR[x_1,\dots,x_n]$ such that the polynomial $q = -g_m^{2\alpha} + \sum_i h_if_i$ is a sum of squares} 
					  \STATE Set $success := $ {\sc True}.
					\ELSE
				      \STATE Increment $\alpha := \alpha + 1$.
				      \IF{$\alpha > \alpha_{\text{max}}$}
				        \RETURN{($S$,{\sc False})}
				      \ENDIF
				    \ENDIF
				  \ENDWHILE
				\ENDFOR
				\RETURN{($S$,{\sc True})}
			\end{algorithmic}
		\end{algorithm}

		
		\section{Equalities and Inequalities}\label{Sec:Inequalities}

One can naturally generalize from real radicals of systems of polynomial equations
to $\sA$-radicals of systems of polynomial equations and inequalities.  
In particular, let $f_1,\dots,f_k,r_1,\dots,r_s\in\bR[x_1,\dots,x_n]$
with
$$\mbox{\scriptsize $I = \langle f_1,\dots,f_k\rangle \hbox{~and~} 
\sA = \{x\in\bR^n~|~r_i(x) \geq 0 \hbox{~for all~}i = 1,\dots,s\}$}.$$
The $\sA$-radical of $I$ is $\sqrt[\sA]{I} = I(\Var_\bR(I)\cap\sA)$.
Algebraically, one can characterize $\sqrt[\sA]{I}$ using sums of squares
\cite{Marshall2008,Stengle1974}:
\begin{equation}\label{eq:Aradical}
\sqrt[\sA]{I} = 
{\scriptsize
\left\{p\in\bR[x]\leftsuchthat\begin{array}{l}
p^{2\alpha} + \displaystyle\sum_{\nu\in\{0,1\}^s} \sigma_\nu\cdot\prod_{j=1}^s r_j^{\nu_j}\in I  \\
\mbox{~~~for some~}\alpha \in \bZ_{> 0}, \\
\hbox{~~~~~~sum of squares~}\sigma_\nu\in\bR[x]\end{array} \right\}\right..}
\end{equation}
Rather than try to locate sample points that satisfy equalities
and inequalities, we will instead reduce to equations by 
introducing ``slack'' variables.  That is,
we consider the ideal
$$J = \langle f_1(x),\dots,f_k(x),r_1(x)-y_1^2,\dots,r_s(x)-y_s^2\rangle.$$
Since $\Var_\bR(I)\cap\sA = \pi(\Var_\bR(J))$ where $\pi(x,y) = x$, we know
\begin{equation}\label{eq:Aradical2}
\sqrt[\sA]{I} = \sqrt[\bR]{J}\cap\bR[x_1,\dots,x_n].
\end{equation}
Thus, we compute $S\subset\overline{\Var_\bR(J)}$
but only perform interpolation on $\pi(S)$.
If $\langle g_1,\dots,g_\ell\rangle = I(\pi(S))\subset\bR[x_1,\dots,x_n]$
and each $g_i\in\sqrt[\bR]{J}$, then $I(\pi(S)) = \sqrt[\sA]{I}$
by \eqref{eq:Aradical2}.

		
		\section{Examples}\label{Sec:Examples}
		
We demonstrate our approach on several examples.
		
		
		\subsection{An illustrative example}\label{Sec:BivariateCubic}
\begin{figure}[!b] \begin{center}
			\includegraphics[scale=0.8,width = 2in]{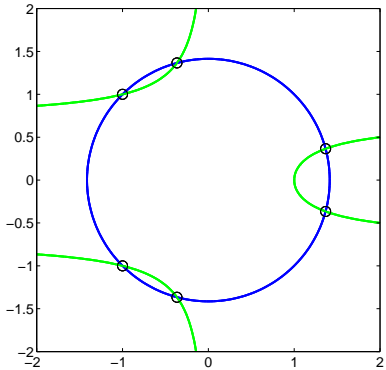}
			\caption{Plot of solutions for $f$ from \S~\ref{Sec:BivariateCubic}}
			\label{fig:example1}
		\end{center}\end{figure}

		To illustrate our approach, we consider the intersection of a circle 
		and a bivariate cubic, namely
		$$f = \{x^2 + y^2 - 2,~2xy^2 - x + 1\}.$$
		The system $f = 0$ has six solutions, all of which are real:
$$\rvar{f} = \{(-1, \pm 1), \, (1.366, \pm 0.366), (-0.366, \pm 1.366)\}$$
		which is shown in Figure~\ref{fig:example1}.
		
In our first test, we simply take $S = \rvar{f}$.  
Since the Hilbert function of $I(S)$ is $1,3,5,6,6,\dots$,
we can show that $I(S)$ is generated by $I(S)_{\leq 3}$.  A basis
for the linear space $I(S)_{\leq 3}$, computed as in \S~\ref{Sec:Interpolate},
is:
$$
		G = \left\{\begin{array}{l} y^3 + x^2y - 2y,~xy^2 - x/2 + 1/2, \\
		x^3 - 3x/2 - 1/2,~x^2 + y^2 - 2\end{array}\right\}.
$$
Using either $f$ or a Gr\"obner basis for $\langle f\rangle$, e.g.,
\begin{equation}\label{eq:GB}
\{x^2 + y^2 - 2,~2xy^2 - x + 1,~2y^4 - 5y^2 - x + 2\},
\end{equation}
every $g\in G$ was found to be in $\rrad{\langle f \rangle}$ 
showing that $S$ is indeed equal to $\rvar{f}$.
		
\subsection*{Incomplete solution set}

Suppose that we take $R = \rvar{f} \cap \{y\geq0\}$.
Since the Hilbert function of $I(R)$ is $1,3,3,\dots$
and $I(R)_{\leq 1} = \{0\}$, we know that $I(R)$
is generated by three quadratics, approximately
\begin{equation}\label{eq:quadrics}
G = \left\{\begin{array}{l}
y^2 - 2.049y - 0.18301x + 0.86603, \\
		xy - 0.18301y - 0.68301x + 1/2\\
		x^2 + 0.18301x + 2.049y - 2.866 \end{array}\right\}.
		\end{equation}
Using $\alpha_{\text{max}} = 5$, we were unable to validate 
that any of the polynomials in $G$
where in $\sqrt[\bR]{\langle f\rangle}$.
In fact, we can show that this is indeed correct
since each polynomial in $G$ is nonzero at each of 
the three points in $\rvar{f}\setminus R$.

		
		\subsection*{Semialgebraic condition}

We now validate that $R = \rvar{f} \cap \{y\geq0\}$ is
the complete solution set for the $\sA$-radical
of $\langle f\rangle$ where $\sA = \{y\geq0\}$.
To that end, we add a slack variable $z$ and consider
the system
$$F = \{x^2 + y^2 - 2, ~2xy^2 - x + 1,~y - z^2\}.$$
As described in \S~\ref{Sec:Inequalities}, we just
need to show that each polynomial in $G$ from \eqref{eq:quadrics}
is contained in $\sqrt[\bR]{\langle F\rangle}$.  
Using either $F$ or a Gr\"obner basis for $\langle F\rangle$,
namely \eqref{eq:GB} together with $y - z^2$, 
we validated that $G\subset\sqrt[\bR]{\langle F\rangle}$
showing that $R$ is indeed equal to $\rvar{f} \cap \{y\geq0\}$,
i.e., $\sqrt[\sA]{\langle f\rangle} = I(R)$.
		
		
		\subsection{Positive-dimensional components}\label{Sec:PosDim}
		
To illustrate the approach on a system such that the real radical
ideal is positive-dimensional, consider the system
		\[
		f = \{xyz, \, z(x^2 + y^2 + z^2 + y), \, y(y + z)\}.
		\]
The set $\Var_\bC(f)$ consists of three lines, 
two of which are complex conjugates of each other that
intersect at the origin and the other is a double line
with respect to $f$, and an isolated point.  
In particular, $\Var_\bR(f)$ is the line
$y = z = 0$ and the isolated point $(0,-1/2,1/2)$.
So, we take 
$$S = \{(x,0,0)~|~x\in\bC\}\cup \{(0,-1/2,1/2)\}\subset\overline{\Var_\bR(f)}.$$

To simplify the real computations later, we first replace
$f$ with a Gr\"obner basis for the radical $\sqrt{\langle f\rangle}$, namely
		\[
		f = \{2yz - y, \, 2y^2 + y, \, xy, \, 4x^2z + 4z^3 + y\}.
		\]
With the isolated solution, sampling 3 points on the 
line is enough
to compute a basis for $I(S)_{\leq 2}$ which generates $I(S)$:
\[
		G = \{ z^2 + y/2, \, yz - y/2, \, y^2 + y/2, \, xz, \, xy, \, y + z \}.
\]
Each element in $G$ was shown to belong to $\rrad{\langle f\rangle}$
with $\alpha \leq 2$.

		
		\subsection{Katsura-5 system}\label{Sec:Katsura}
		
As an illustration of our approach on a problem which
was solved using the semidefinite characterization
of the real radical in \cite{LLR08}, we consider
the Katsura-5 system as in \mbox{\cite[Ex.~5.4]{LLR08}}.
The system consists of a linear, say $f_1$, 
and five quadratics, say $f_2,\dots,f_6$, in 
six variables, namely
$$f = 
\hbox{\scriptsize $\left\{\begin{array}{c}
x_1 + 2(x_2 + x_3 + x_4 + x_5 + x_6) - 1, \\
x_1^2 + 2(x_2^2 + x_3^2 + x_4^2 + x_5^2 + x_6^2) - x_1, \\
2(x_1x_2 + x_2x_3 + x_3x_4) + x_4x_5 + x_5x_6 -x_2, \\
x_2^2 + 2(x_1x_3 + x_2x_4 + x_3x_5 + x_4x_6) - x_3, \\
2(x_1x_4 + x_2x_3 + x_2x_5 + x_3x_6) - x_4, \\
x_3^2 + 2(x_1x_4 + x_1x_5 + x_1x_6) - x_5
\end{array}\right\}$}.$$
The set $\Var_\bC(f)$ consists of $32$ points,
$12$ of which lie in $\bR^6$.  
The set of real solutions, say $S$, is readily computed 
using homotopy continuation.  

The Hilbert function is $1,6,12,12,\dots$ 
with $I(S)$ being generated by $I(S)_{\leq 2}$.
In particular, $I(S)_{\leq 2}$ is a linear space
spanned by the linear $f_1$ and $15$ quadratics\footnote{Available
at \url{www.nd.edu/~aliddel1/validate-reals}.}.

Trivially, $f_1\in\sqrt[\bR]{\langle f_1,\dots,f_6\rangle}$
and the quadratics are shown to be in the real radical using
$\alpha\leq 2$.  This computation validates that $\Var_\bR(f)$ consists
of $12$ points.  Moreover, this data matches
that displayed in \cite[Table~4]{LLR08}.


\subsection{Seiler system}\label{Sec:Seiler}

As an illustration of our approach on a problem 
considered in \cite[Ex.~5]{ma2014certificate}, 
namely the Seiler system \cite{Seiler}
$$f = \left\{\begin{array}{c} x_3^2 + x_2 x_3 - x_1^2, \\
x_1x_3 + x_1x_2 - x_3,  \\
x_2x_3 +x_2^2 + x_1^2 -x_1\end{array}\right\}.$$
This system does not have a Pommaret basis with respect
to the total degree ordering defined by 
$x_1 < x_2 < x_3$ \cite{Seiler}.  Thus,~\cite{ma2014certificate}
uses a change of coordinates to overcome this.

Even though $f$ consists of $3$ polynomials in $3$ variables,
$\Var_\bC(f)$ is actually a curve.  In particular, 
$I = \langle f\rangle$ is a one-dimensional prime ideal, i.e.,
$I = \sqrt{I}$ and $\Var_\bC(I)$ is an irreducible curve.
Hence, we know that $I = \sqrt[\bR]{I}$ if we can compute
a real point $x\in\Var_\bR(I)$ which is smooth with respect
to $f$, i.e., the rank of $Jf(x)$ is $2$.  

To that end, we utilize a gradient descent homotopy \cite{GH15}.
We took $y = (1,-3/2,3/4)$ and considered the homotopy
$$\hbox{\small $H(x,\lambda,t) = 
\left[\begin{array}{c}
f(x) - t\cdot f(y) \\
\lambda_0 (x-y) + \lambda_1 \nabla f_1(x) + 
\lambda_2 \nabla f_2(x) + \lambda_3 \nabla f_3(x) \end{array}\right]$}$$
where $\lambda\in\bP^3$.  Starting at
$x = y$ and \mbox{$\lambda = [1,0,0,0]\in\bP^3$}
when $t = 1$, we obtain a point, which is
approximately $(0.7009,-0.2504,-0.5868)$, that lies on
$\Var_\bR(f)$ and is indeed a smooth point on $\Var_\bC(f)$.
Hence, the isosingular set of this point with respect to $f$
is $\Var_\bC(f)$ showing that $I = \sqrt[\bR]{I}$.
		
		
		\subsection{An energy landscape}\label{Sec:Phi4N3}
		
\begin{figure}[!t]
			\centering
			\begin{tikzpicture}[<->,>=stealth', auto, node distance=1.5cm,
			thick,main node/.style={circle,draw,font=\sffamily\bfseries},scale=0.5]
			\clip(-2.5,-7.5) rectangle (7.5,2.5);
			\node[main node] (1) {$1$};
			\node[main node] (2) [right of=1] {$2$};
			\node[main node] (3) [right of=2] {$3$};
			\node[main node] (4) [below of=1] {$4$};
			\node[main node] (5) [right of=4] {$5$};
			\node[main node] (6) [right of=5] {$6$};
			\node[main node] (7) [below of=4] {$7$};
			\node[main node] (8) [right of=7] {$8$};
			\node[main node] (9) [right of=8] {$9$};
			\path[every node/.style={font=\sffamily\small}]
			(1) edge node [right] {} (2)
			(2) edge node [right] {} (3)
			(1) edge node [right] {} (4)
			(2) edge node [right] {} (5)
			(3) edge node [right] {} (6)
			(4) edge node [right] {} (5)
			(5) edge node [right] {} (6)
			(4) edge node [right] {} (7)
			(5) edge node [right] {} (8)
			(6) edge node [right] {} (9)
			(7) edge node [right] {} (8)
			(8) edge node [right] {} (9)
			(1) edge[dashed,bend left=150] node [left] {} (3)
			(1) edge[dashed,bend left=150]  node [left] {} (7)
			(4) edge[dashed,bend left=150] node [left] {} (6)
			(2) edge[dashed,bend left=150]  node [left] {} (8)
			(7) edge[dashed,bend left=150] node [left] {} (9)
			(3) edge[dashed,bend left=150]  node [left] {} (9) ;
			\end{tikzpicture}
			\caption{Nearest-neighbor coupling for a $3\times3$ grid of nodes.}
			\label{fig:grid} 
		\end{figure}
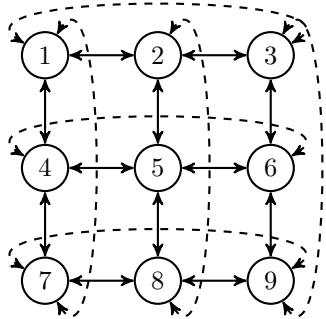

Our final example aims to compute the real critical points of
the energy landscape of the two-dimensional nearest-neigh\-bor 
$\phi^4$ model on a $3\times3$ grid as in \cite{Franzosi,EnergyLandscape}.
We label the nodes $1,\dots,9$ with Figure~\ref{fig:grid} showing
the coupling between the nodes.  Let $N(i)$ denote the four 
nearest neighbors of node $i$, e.g., $N(1) = \{2,3,4,7\}$.
After selecting various parameters for this model, we consider 
the potential energy
$$\mbox{$\displaystyle V(x) = \sum_{i=1}^9\left[\frac{1}{40} x_i^4 - x_i^2 + \frac{1}{4}\sum_{j\in N(i)} (x_i-x_j)^2\right]$}.$$
The system defining the critical points is $f = \nabla V$ so that
$$f_i = \frac{1}{10}x_i^3 - 2x_i + \sum_{j\in N(i)}(x_i - x_j)$$

The system $f$ is a Gr\"obner basis and the set $\Var_\bC(f)$ 
consists of $3^9 = \mbox{19,683}$ points.  
However, when searching for real stationary points, one
only obtains $3$ points, namely 
$$S = \{(0, 0, 0, 0, 0, 0, 0, 0, 0),~\pm (w, w, w, w, w, w, w, w, w)\}$$
where $w = \sqrt{20}\approx 4.4721$.  Hence,
$I(S)$ is generated by
		\[
		G = \{ x_1(x_1^2 - 20), x_2 - x_1, \dots, x_9 - x_1\}.
		\]
All nine basis elements were found to be in $\rrad{\langle f \rangle}$
with $\alpha = 1,2,\dots,2$, respectively.
Therefore, $S = \Var_\bR(f)$, i.e., the energy landscape
$V$ has exactly three real critical points.
			
		\section{Conclusion} \label{Sec:Conclusion}

By combining numerical algebraic geometry with sums of squares
programming, we have produced a method for certifying 
that a set of polynomials generate the real radical.
The set of polynomials arises from the generators of a 
set $S$ which is contained in the Zariski closure of the set of real solutions.
As first considered in \cite{CP2015}, combining numerical algebraic
geometry and semidefinite programming can improve
the efficiency of computations and produce new approaches, in particular
for computing and analyzing the set of real solutions of a system of 
polynomial equations. 

		\bibliographystyle{plain}
		\vskip 0.2in
		\small
		\bibliography{citations}
	\end{document}